\documentclass[11pt,a4paper]{article}
\usepackage{a4}
\usepackage{graphics}

\usepackage{amssymb,epsfig}
\usepackage{amsmath}
\usepackage{amsthm}
\usepackage{color}
\usepackage{enumerate}
\usepackage{psfrag}

\def\div{{\rm div \;}}

\def\N{\mathbb{N}}
\def\R{\mathbb{R}}

\def\E{\mathbb{E}} 
\def\P{\mathbb{P}} 

\newtheorem{theorem}{Theorem}
\newtheorem{proposition}[theorem]{Proposition}

\newtheorem{definition}[theorem]{Definition}
\newtheorem{lemma}[theorem]{Lemma}
\newtheorem{remark}[theorem]{Remark}

\begin{document}
\title{Accelerated dynamics: Mathematical foundations and algorithmic improvements}
\author{Tony
  Leli\`evre\thanks{Universit\'e Paris-Est, CERMICS (ENPC), INRIA, 6-8 Avenue Blaise
Pascal, F-77455 Marne-la-Vall\'ee. The work of
    T. Leli\`evre is supported by the European Research Council under
    the European Union's Seventh Framework Programme (FP/2007-2013) /
    ERC Grant Agreement number 614492. T. Leli\`evre would like to
    thank very instructive discussions with D. Perez and A.F. Voter on
  accelerated dynamics, as well as his co-authors D.~Aristoff, A.~Binder,
  C.~Le~Bris, M.~Luskin, F.~Nier, D.~Perez and G.~Simpson.}}
\maketitle
\begin{abstract}
We present a review of recent works on the mathematical analysis of
algorithms which have been proposed by A.F. Voter and co-workers in
the late nineties in order to efficiently generate long trajectories
of metastable processes. These techniques have been successfully
applied in many contexts, in particular in the field of materials
science. The mathematical analysis we propose relies on the notion of
quasi stationary distribution.
\end{abstract} 
\section{Introduction}
\label{intro}

This article is a review of recent works whose aim is to lay the
mathematical foundations of algorithms used in computational
statistical physics, namely {\em accelerated dynamics techniques}
introduced by A.F.~Voter and co-workers in the late nineties. These methods have been
proposed in order to efficiently sample trajectories in the context of
molecular dynamics.

Molecular dynamics is used in various application fields (biology,
chemistry, materials science) in order to simulate the evolution of a
molecular system, namely interacting particles representing atoms or
group of atoms. The typical dynamics one should have in mind is the
Langevin dynamics:
\begin{equation}\label{eq:Langevin}
\left\{
\begin{aligned}
dq_t & = M^{-1} p_t \, dt \\
dp_t &= -\nabla V(q_t) \, dt - \gamma M^{-1} p_t \, dt + \sqrt{2
  \gamma \beta^{-1}} dW_t
\end{aligned}
\right.
\end{equation}
where $(q_t,p_t)$ denotes the positions and momenta of the particles
at time $t \ge 0$,
$M$~is the mass tensor,
$V$ is the potential function which, to a given set of positions $q$,
associates its energy $V(q)$, $\gamma >0$ is a friction parameter,
$\beta = (k_B T)^{-1}$ is proportional to the inverse temperature and $W_t$ is a
standard Brownian motion. In the following, we assume that $q_t \in
\R^d$ where $d$ is typically very large (say $3$ times the number of
particles) but generalizations to dynamics on manifolds (systems with
constraints) are straightforward.
When $\gamma=0$, the Langevin dynamics is nothing but the Hamiltonian dynamics. The terms involving $\gamma$ model the fact that the system is at a given
temperature. Indeed, under loose
assumptions on $V$, this dynamics is ergodic with respect to the
canonical measure (NVT ensemble): for any test function $\varphi:\R^d
\times \R^d \to \R$,
$$\lim_{T \to \infty} \frac{1}{T} \int_0^T \varphi (q_t,p_t) \, dt =
Z^{-1} \int \varphi(q,p) \exp (-\beta (p^T M^{-1} p /2 + V(q))) \, dp
dq$$
where $Z= \int \exp (-\beta (p^T M^{-1} p /2 + V(q))) \, dp
dq< \infty$. In the following, we will also consider the overdamped
Langevin dynamics which is obtained from~\eqref{eq:Langevin} in the
limit $\gamma \to \infty$ or $M \to 0$ (see for example~\cite[Section 2.2.4]{lelievre-rousset-stoltz-book-10}):
\begin{equation}\label{eq:overdamped_Lang}
dX_t=-\nabla V(X_t) \, dt + \sqrt{2 \beta^{-1}} dW_t
\end{equation}
where $X_t \in \R^d$ denotes the position of the particles. Again, under loose
assumptions on $V$, this dynamics is ergodic with respect to the
canonical measure $\tilde{Z}^{-1} \exp(-\beta V(x)) \, dx$ where $\tilde{Z}= \int \exp (-\beta V(x)) \, dx$. The aim of
molecular simulations is to compute macroscopic properties from the
models~\eqref{eq:Langevin} or~\eqref{eq:overdamped_Lang} at the
atomistic level ($V$ being the main modelling ingredient). In this
article, we are particularly interested in so-called dynamical
quantities, namely macroscopic observables which depend on the path
$(q_t,p_t)_{t \ge 0}$ or $(X_t)_{t \ge 0}$. For example, one would
like to sample the paths which go from one region of the phase
space to another one, in order to compute the typical time to observe
such transitions or the intermediate states along the transition
path.

The numerical difficulty associated with such computations is that the
timescale at the microscopic level is much smaller than the timescale
at the macroscopic level. More precisely, the timestep required to
obtain a stable discretization of the above dynamics is much smaller
than the timescale associated with the macroscopic observables  of interest. In
other words, one has to simulate very long trajectories of a high
dimensional stochastic dynamics. In practice, for many applications, a
naive discretization of the dynamics is not sufficient to reach the
timescales of interest, since it would require up to typically $10^{15}$
iterations (the timescale at the atomistic level - bond vibration - is indeed of the order of
$10^{-15}\, {\rm s}$, while transitions between metastable regions
occur may over timescales ranging from microseconds to seconds).

The idea of accelerated dynamics is to take benefit from this
timescale discrepancy in order to simulate more efficiently paths over
very long times. Indeed, the typical trajectories
of~\eqref{eq:Langevin} or~\eqref{eq:overdamped_Lang} are metastable:
this means that the trajectory remains trapped for very long times in
some region of the phase space, before hopping to another region where
it again remains trapped. These regions are called metastable
states. Metastability originates from energetic barriers (the path to
leave the state requires to climb above a saddle point of the
potential energy $V$) or from entropic barriers (the path to
leave the state goes through a narrow corridor, due to some steric
constraints in the system for example), or more generally from a
combination of energetic and entropic effects. The bottom line is thus that behind the continuous state space
dynamics~\eqref{eq:Langevin} or~\eqref{eq:overdamped_Lang}, there is a
discrete state space jump process (encoding the jumps from metastable
states to metastable states). Actually, discrete state space Markov
dynamics are also very much used in molecular dynamics: there are
called kinetic Monte Carlo or Markov state models, see for
example~\cite{voter-05}. And continuous state space models are typically used
in order to parametrize these Markovian models (namely to compute the
jump rates between metastable states) using for example Arrhenius
(or Eyring-Kramers) formulas. The accelerated dynamics of A.F.~Voter
follow a different path: the principle is to use the underlying jump
process in order to accelerate the sampling of the original dynamics,
in the spirit of a predictor-corrector schemes. These are thus
numerical methods to efficiently generate the underlying jump
process among metastable states.

In the following, we will assume that we are given a mapping 
\begin{equation}
\mathcal{S}: \R^d \to \N
\end{equation}
which to a given set of positions $x \in \R^d$ associates
$\mathcal{S}(x)$, the label of the state in which $x$ lies. One
should think of the states
$$\mathcal{S}^{-1}(\{n\})=\{x \in \R^d, \mathcal{S}(x)=n\} \text{ for } n
\in \N$$
as the metastable states mentioned above. This mapping thus defines a
partition of the state space.
Let us
make two comments on this mapping. First, an important message from
the mathematical analysis we present below is that whatever the mapping
$\mathcal{S}$, the accelerated dynamics algorithms are consistent: they give
the correct result in some limiting regime. For example for the
parallel replica method, in the limit when the correlation time - a
numerical parameter introduced below - goes to infinity, the generated
dynamics are statistically correct. In
particular if some of the states happen not to be metastable, or if
for one specific realization, the stochastic process does not remain trapped in
one of this state (because, for example, it enters the state with a
too large velocity), the algorithms are still consistent. Second, as
will become clear below, the numbering of the states do not need to be
known a priori: the states are numbered as the simulation goes, when
they are successively discovered by the stochastic process.

Three algorithms have been proposed by A.F. Voter and co-workers. The idea is that if
the stochastic process remains trapped for a very long time in a given state
$S=\mathcal{S}^{-1}(\{n\})$ (for some given $n$), then there are ways to generate the exit
event from this state much more efficiently than by running the
original dynamics until the exit time. The exit event is fully
characterized by two random variables: the exit time and the exit
point from $S$, which are defined by considering
the first hitting time and point on the boundary $\partial S$. Let us
roughly describe the ideas behind the three algorithms.

 In the
{\em Parallel Replica method}~\cite{voter-98}, the principle is to simulate
in parallel many trajectories following the original
dynamics~\eqref{eq:Langevin} or~\eqref{eq:overdamped_Lang}, to
consider the first exit event among the replicas, and to
generate from this first exit event a consistent exit time and exit
point. The gain is thus obtained in terms of wall clock time. This
algorithm can be seen as a way to parallelize a computation in time,
which is not an easy problem in general due to the sequential nature
of time evolutions. 

In the {\em hyperdynamics}~\cite{voter-97}, the idea is to modify
the potential $V$ within the state~$S$ in order to accelerate the exit
from the state for the original dynamics~\eqref{eq:Langevin} or~\eqref{eq:overdamped_Lang}. Again, using an appropriate time rescaling, it is
possible to generate from the observed exit event on the biased
potential an exit event which is consistent with what would have been
observed on the original unbiased potential $V$. 

The {\em Temperature Accelerated Dynamics}
(TAD)~\cite{sorensen-voter-00} consists in considering the original
dynamics~\eqref{eq:Langevin} or~\eqref{eq:overdamped_Lang} at a higher
temperature than the original one. The idea is then that under
appropriate assumptions, there is a way to infer from the exit events
observed at high temperature the exit event which would have been
observed at the original lower temperature.

The ultimate aim of these three techniques is thus to generate efficiently
the so-called state-to-state dynamics $({\mathcal S}(q_t))_{t \ge 0}$
(for~\eqref{eq:Langevin}) or $({\mathcal S}(X_t))_{t \ge 0}$
(for~\eqref{eq:overdamped_Lang}), with the correct statistical
properties. Let us emphasize that the objective is to get the correct
law on the paths (in order to compute dynamical quantities), not only
on the time marginals for example.

A crucial mathematical tool to understand these techniques is the
Quasi-Stationary Distribution (QSD) introduced in
Section~\ref{sec:QSD}. We will then describe the mathematical results
which have been obtained so far on the three algorithms: Parallel
Replica in Section~\ref{sec:ParRep}, hyperdynamics in
Section~\ref{sec:Hyper} and Temperature Accelerated Dynamics in
Section~\ref{sec:TAD}. A few concluding remarks are provided in Section~\ref{sec:conc}.

\section{The Quasi-Stationary Distribution and the decorrelation
  step}\label{sec:QSD}

Let us consider a fixed state $S=\mathcal{S}^{-1}(\{n\})$ and let us
focus for simplicity on the overdamped Langevin
dynamics~\eqref{eq:overdamped_Lang}. We assume that $S$ is a bounded
regular subset of $\R^d$. Let us consider the first exit
time from $S$:
$$T_S=\inf \{t \ge 0, \, X_t \not\in S \}.$$

\subsection{The QSD}

Let us start with the definition of the QSD.
\begin{definition}
A probability measure $\nu$ with support in $S$ is called a QSD for
the Markov process $(X_t)_{t \ge 0}$ if and only if 
$$\forall t > 0, \, \forall A \subset S, \nu(A)=\frac{\int_S \P(X_t^x
  \in A, t < T^x_S) \, \nu(dx)}{\int_S \P(t < T^x_S) \, \nu(dx)}.$$
\end{definition}
In other words, $\nu$ is a QSD if, when $X_0$ is distributed according
to $\nu$, the law of $X_t$ conditionally to the fact that $(X_s)_{0
  \le s \le t}$ remains in the state $S$ is still $\nu$, for all
positive~$t$.

The QSD satisfies three properties which will be crucial in the
following. We refer for example to~\cite{le-bris-lelievre-luskin-perez-12} for a proof of these results and to~\cite{collet-martinez-san-martin-13} for more
general results on QSDs.

\begin{proposition}\label{prop:QSD_1}
Let $(X_t)_{t \ge 0}$ follow the dynamics~\eqref{eq:overdamped_Lang}
with an initial condition $X_0 \in S$. Then, there exists a probability
distribution $\nu$ with support in $S$ such that
\begin{equation}\label{eq:QSD_1}
\lim_{t \to \infty} {\mathcal L}(X_t | T_S > t) = \nu.
\end{equation}
The distribution $\nu$ is the QSD associated with $S$.
\end{proposition}
A consequence of this proposition is the existence and uniqueness of
the QSD. The QSD can thus be seen as the longtime limit of the process
conditioned to stay in the well. This proposition can be useful to
understand what is a metastable state. A metastable state is a state such that the typical exit time is
much larger than the local equilibration time, namely the time to
observe the convergence to the QSD in~\eqref{eq:QSD_1}.

Let us now give a second property of the QSD.
\begin{proposition}\label{prop:QSD_2}
Let $L=-\nabla V \cdot \nabla + \beta^{-1} \Delta$ be the
infinitesimal generator of $(X_t)_{t \ge 0}$ (satisfying~\eqref{eq:overdamped_Lang}). Let us consider the
first eigenvalue and eigenfunction associated with the adjoint operator
$L^*=\div ( \nabla V + \beta^{-1} \nabla)$ with homogeneous Dirichlet
boundary condition:
\begin{equation}\label{eq:QSD_2}
\left\{
\begin{aligned}
L^* u_1 &= - \lambda_1 u_1 \text{ on $S$},\\
u_1 &=0 \text{ on $\partial S$}.
\end{aligned}
\right.
\end{equation}
The QSD $\nu$ associated with $S$ satisfies:
$$d\nu = \frac{u_1(x) \, dx}{\int_S u_1(x) \, dx}$$
where $dx$ denotes the Lebesgue measure on $S$.
\end{proposition}
The QSD thus has a density with respect to the Lebesgue measure, which
is nothing but the ground state of the Fokker-Planck operator $L^*$
associated with the dynamics with absorbing boundary conditions.

Finally, the last property of the QSD concerns the exit event, when
$X_0$ is distributed according to $\nu$.
\begin{proposition}\label{prop:QSD_3}
Let us assume that $X_0$ is distributed according to the QSD $\nu$ in $S$. Then the
law of the couple $(T_S,X_{T_S})$ (namely the first exit time and the
first exit point) is fully characterized by the following properties:
\begin{itemize}
\item $T_S$ is exponentially distributed with parameter $\lambda_1$
  (defined in Equation~\eqref{eq:QSD_2} above)
\item $T_S$ is independent of $X_{T_S}$
\item The law of $X_{T_S}$ is given by: for any bounded measurable
  function $\varphi: \partial S \to \R$,
\begin{equation}\label{eq:QSD_3}
\E^\nu(\varphi(X_{T_S})) = - \frac{\int_{\partial S}
  \varphi \, \partial_n u_1 \, d \sigma}{\beta \lambda_1 \int_{S}
  u_1(x) \, dx}
\end{equation}
where $\sigma$ denotes the Lebesgue measure on $\partial S$ induced by
the Lebesgue measure in $\R^d$ and the Euclidean scalar product, and
$\partial_n u_1=\nabla u_1 \cdot n$ denotes the outward normal
derivative of $u_1$ on $\partial S$. 
\end{itemize}
\end{proposition}
This Proposition explains the interest of the QSD. Indeed, if the
process is distributed according to the QSD in $S$ (namely, from
Proposition~\ref{prop:QSD_1}, if it remained for a sufficiently long
time in $S$), then the exit event from the state $S$ is Markovian, in
terms of state-to-state dynamics. This is reminiscent of what is
assumed to build kinetic Monte Carlo models (see~\cite{voter-05}).


\begin{remark}
The existence of the QSD and the convergence of the
constrained process towards the QSD for the Langevin
process~\eqref{eq:Langevin} requires extra work compared to the
overdamped Langevin process~\eqref{eq:overdamped_Lang}. For results in that direction, we refer to the recent
manuscript~\cite{nier-14}.
\end{remark}

\subsection{The decorrelation step}\label{sec:decor}

The accelerated dynamics algorithms will be based on the assumption
that the process remained sufficiently long in the state $S$ so that
one can consider it is distributed according to the QSD $\nu$. Then,
using this assumption, various techniques are used in order to
efficiently generate the exit event from $S$, starting from the QSD
(see the next sections).

A natural preliminary question is therefore: how to assess in practice that the
limit has been reached in~\eqref{eq:QSD_1} ? This is done in the
so-called decorrelation step which consists in waiting for a given
time $\tau_{corr}$ (a so-called decorrelation time) before assuming
that the local equilibrium $\nu$ has been reached. This correlation
time can be state dependent, and is typically supposed to be known {\em
  a priori}. 

From a mathematical viewpoint, $\tau_{corr}$ should be
chosen sufficiently large so that the distance between the law of
$X_{\tau_{corr}}$ conditioned to $T_S \ge \tau_{corr}$ and the QSD $\nu$
is small. In~\cite{le-bris-lelievre-luskin-perez-12}, we prove the following:
\begin{proposition}\label{prop:CV_QSD}
Let $(X_t)_{t \ge 0}$ satisfies~\eqref{eq:overdamped_Lang} with $X_0
\in S$. Let us consider $-\lambda_2 < -\lambda_1 < 0$  the first  two
eigenvalues of the operator $L^*$ on $S$ with homogeneous Dirichlet boundary
conditions on $\partial S$ (see Proposition~\ref{prop:QSD_2} for the
definition of $L^*$). Then, there
exists a constant $C>0$ which depends on the law of $X_0$, such that,
for all $t \ge \frac{C}{\lambda_2  - \lambda_1}$,
$$\sup_{f,\, \|f\|_{L^\infty} \le 1} \left| \E(f(T_S-t,X_{T_S}) | T_S
  \ge t) -
  \E^\nu(f(T_S,X_{T_S})) \right| \le C \exp (-(\lambda_2-\lambda_1) t ).$$
\end{proposition}
In other words, the total variation norm between the law of
$(T_S-t,X_{T_S})$ conditioned to $T_S \ge t$ (for any initial
condition $X_0 \in S$), and the law of
$(T_S,X_{T_S})$ when $X_0$ is distributed according to $\nu$,
decreases exponentially fast with rate $\lambda_2-\lambda_1$. This
means that $\tau_{corr}$ should be chosen of the order
$1/(\lambda_2-\lambda_1)$. Of course, this is not a very practical
result since computing these eigenvalues is in general
impossible. From a theoretical viewpoint, this result tells us that
the local equilibration time is of the order
$1/(\lambda_2-\lambda_1)$, so that, the state $S$ will be metastable
if this time is much smaller than the exit time (which is typically of
the order of $1/\lambda_1$, see Proposition~\ref{prop:QSD_2}). 

Let us mention the recent work~\cite{binder-simpson-lelievre-14} where we propose a
numerical method to
approximate the time to reach the QSD using a Fleming-Viot particle process together
with stationarity statistical tests. The interest of the approach is
demonstrated on toy examples (including the 7 Lennard-Jones cluster),
but it remains to test this technique on higher dimensional problems.

From now on, we assume that the decorrelation step has been
successful, and we look for efficient techniques to generate the exit
event (namely a sample of the random variables $(T_S,X_{T_S})$). Let
us describe successively the three algorithms which has been proposed
by A.F. Voter and co-workers.

\section{The Parallel Replica method}\label{sec:ParRep}

\begin{figure}[htbp]\centerline{\includegraphics[width=7cm]{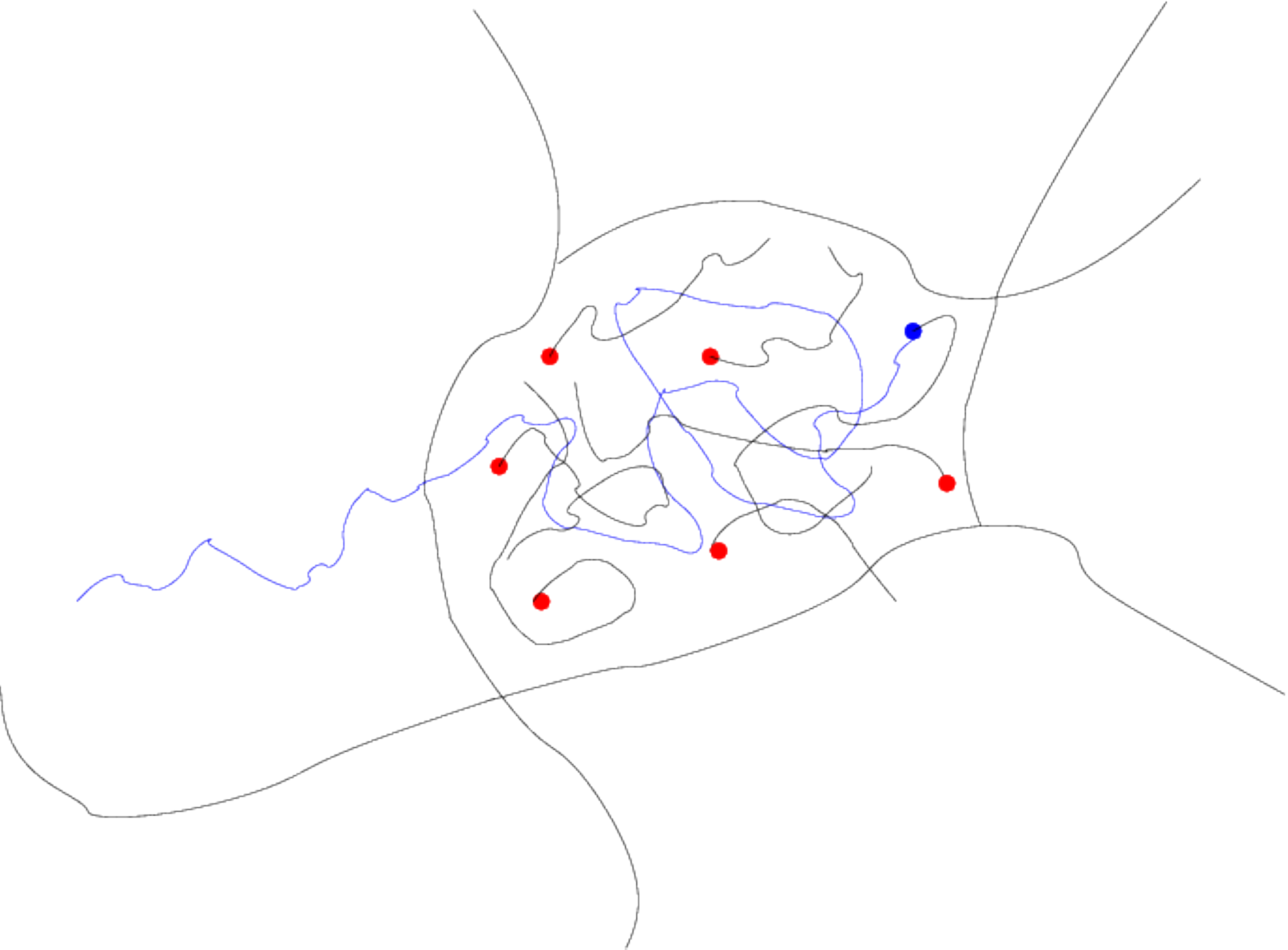}}
\caption{The Parallel Replica method: many exit events are simulated
  in parallel, all starting from the QSD in $S$. The blue trajectory
  represents the reference walker which stays sufficiently long within
$S$ so that we can assume the blue point is distributed according to
the QSD. The red points represent i.i.d. initial conditions
distributed according to the QSD. The black trajectories are simulated
in parallel.}
\end{figure}\label{fig:ParRep}

Let us assume that we are given an initial condition $X_0 \in S$ such
that $X_0$ is distributed according to the QSD in $S$. Let us assume
that we are given a computer with many CPUs in parallel.
The idea of the parallel replica method is to distribute $N$
independent initial conditions $(X_0^i)_{1 \le i \le N}$ in $S$
according to the QSD $\nu$, to let them evolve according
to~\eqref{eq:overdamped_Lang} driven by independent motions (so that
the replicas remain independent) and
then to consider the first exit event among the replicas:
\begin{equation}\label{eq:I0}
I_0=\arg\min_{i \in \{1, \ldots,N\}} T_S^i \text{ where }
T_S^i=\inf\{t \ge 0,\, X^i_t \not\in S\}.
\end{equation}
The integer $I_0 \in \{1, \ldots ,N\}$ is the index of the first replica which exits $S$, and $\min(T^1_S, \ldots , T^N_S)=T^{I_0}_S$. The effective
exit time is set as $N$ times the first exit time, and the effective
exit point is nothing but the exit point for the first exit event. See
Figure~\ref{fig:ParRep} for a schematic illustration of the method.

The consistency of the method is a corollary of
Proposition~\ref{prop:QSD_3}. Indeed using the fact that, starting
from the QSD, the exit time is exponentially distributed and
independent of the exit point, we easily obtain that
\begin{equation}\label{eq:Nminexpo}
N T^{I_0}_S = N \min(T^1_S, \ldots, T^N_S) \stackrel{\mathcal
  L}{=}T^1_S
\end{equation}
which means that the effective exit time has the correct law and
$$X^{I_0}_{T^{I_0}_S}\stackrel{\mathcal L}{=}X^1_{T^1_S}$$
which means that the first exit point of the replica $I_0$ (the first
one to exit among $N$) has the same law as the first exit point of any
of them. Moreover $N T^{I_0}_S=N \min(T^1_S, \ldots , T^N_S)$ and
$X^{I_0}_{T^{I_0}_S}$ are independent, so that we have proven the
following Lemma.
\begin{lemma}
Let $I_0$ be the index of the first replica exiting $S$, defined
by~\eqref{eq:I0}. Then we have the equality in law:
$$\left(N T^{I_0}_S,X^{I_0}_{T^{I_0}_S}\right) \stackrel{\mathcal L}{=} (T^1_S,X^1_{T^1_S}).$$
\end{lemma}
This Lemma shows that the parallel replica is exact: the law of the
effective exit time and exit point is exactly the law of the exit time
and exit point which would have been observed for only one replica.

Let us make a few remarks on this algorithm. First, the full algorithm
actually iterates three steps:
\begin{itemize}
\item The {\em decorrelation step} (see Section~\ref{sec:decor}), where a
  reference walker is run following the
  dynamics~\eqref{eq:overdamped_Lang} until it remains trapped for a
  sufficiently long time in one of the sets ${\mathcal
    S}^{-1}(\{n\})$, so that is can be assumed to be distributed
  according to the QSD $\nu$ associated with this set. During this step,
  the algorithm thus consists simply in integrating the original
  dynamics. No error is introduced and there is no computational gain.
\item The {\em dephasing step} which is a preparation step during which $N$
  independent initial conditions distributed according to $\nu$ are
  sampled, each one on a different CPU. This is done in
  parallel. During this step, the simulation clock is stopped. This
  step is thus pure overhead. This step requires
  appropriate algorithms to sample the QSD such as rejection algorithm
  or Fleming-Viot particle systems (see~\cite{le-bris-lelievre-luskin-perez-12}). For
  example, the rejection algorithm consists in running independently
  walkers following the dynamics~\eqref{eq:overdamped_Lang} (starting
  from a point within $S$) and to
  consider the final point of the trajectory conditionally to the fact
  that the walker remains in the state, for a sufficiently long
  trajectory (typically for the time $\tau_{corr}$ introduced in Section~\ref{sec:decor}).
\item The {\em parallel step}, just described above, which consists in
  running the $N$ replicas independently in parallel, and in waiting
  for the first exit event among the $N$ replicas. The simulation
  clock is then updated by adding the effective exit time
  $NT^{I_0}_S$. The exit point $X^{I_0}_{T^{I_0}_S}$ is used as the
  initial condition of the reference walker for the next
  decorrelation step. The computational gain of the whole algorithm
  comes from this step which divides the wall clock time to sample the
  exit event by the number of
  replicas $N$.
\end{itemize}
In practice, if the rejection algorithm is used in the dephasing step,
one actually does not need to wait for the $N$ replicas
to be dephased to proceed to the parallel step, see~\cite{voter-98,le-bris-lelievre-luskin-perez-12,binder-simpson-lelievre-14}.

In view of the above discussions, the errors introduced in the
algorithm have two origins. First, in the decorrelation step, $\tau_{corr}$
should be chosen sufficiently large so that at the end of the
decorrelation step, the reference walker is indeed distributed
according to a probability law sufficiently close to the QSD. The
convergence result of Proposition~\ref{prop:CV_QSD} shows that the
error is of the order
$O(\exp(-\tau_{corr}/(\lambda_2-\lambda_1)))$. Second, in the
dephasing step, the sampling algorithm of the QSD should be
sufficiently precise in order to obtain i.i.d. samples distributed
according to $\nu$. For the rejection algorithm, independence is
ensured, and the accuracy is again related to the convergence result
of Proposition~\ref{prop:CV_QSD}. For Fleming-Viot particle process,
some correlations are introduced among the replicas, and it is an open
problem to evaluate the error introduced by these correlations. As
already mentioned above, in~\cite{binder-simpson-lelievre-14}, we recently proposed an algorithm
to compute on-the-fly a good correlation time, while sampling the
QSD, using a Fleming-Viot particle process.

The parallel replica is thus a very versatile algorithm. In particular
it applies to both energetic and entropic barriers. The only errors
introduced in the algorithm are related to the rate of convergence of the
conditioned process to the QSD. The algorithm will be all the more
efficient than the convergence time to the QSD is small compared to
the exit time (namely the states are metastable): in this case, the speed-up in terms of wall clock time is linear as a function
of $N$. We refer to~\cite[Section 5.1]{binder-simpson-lelievre-14} for
a discussion of the parallel efficiency of the algorithm.

Let us finally mention the recent work~\cite{aristoff-lelievre-simpson-14} where we propose an
extension of the Parallel Replica algorithm to Markov chains (namely discrete-in-time
stochastic processes). This is indeed a relevant question since in
practice, the continous-in-time dynamics (such as~\eqref{eq:Langevin} or~\eqref{eq:overdamped_Lang}) are approximated by
discrete-in-time Markov Chains using appropriate time discretization schemes.
 The algorithm has to be slightly adapted since,
starting from the QSD (which is still perfectly well defined in this
context), the exit time is not exponentially distributed but has a geometric
law. Therefore, the formula~\eqref{eq:Nminexpo} does not hold and is replaced by
the following: for $T^1, \ldots, T^N$ $N$ i.i.d. random variables
geometrically distributed,
$$N ( \min(T^1,
\ldots, T^N)-1 ) + \min(i \in \{1, \ldots ,N \} , \, T^i =  \min(T^1,
\ldots, T^N) ) \stackrel{\mathcal
  L}{=}T^1.$$
This yields a natural adaptation of the original Parallel Replica
algorithm to Markov chains.

\section{The hyperdynamics}\label{sec:Hyper}

\begin{figure}[htbp]
\centerline{\includegraphics[width=7cm]{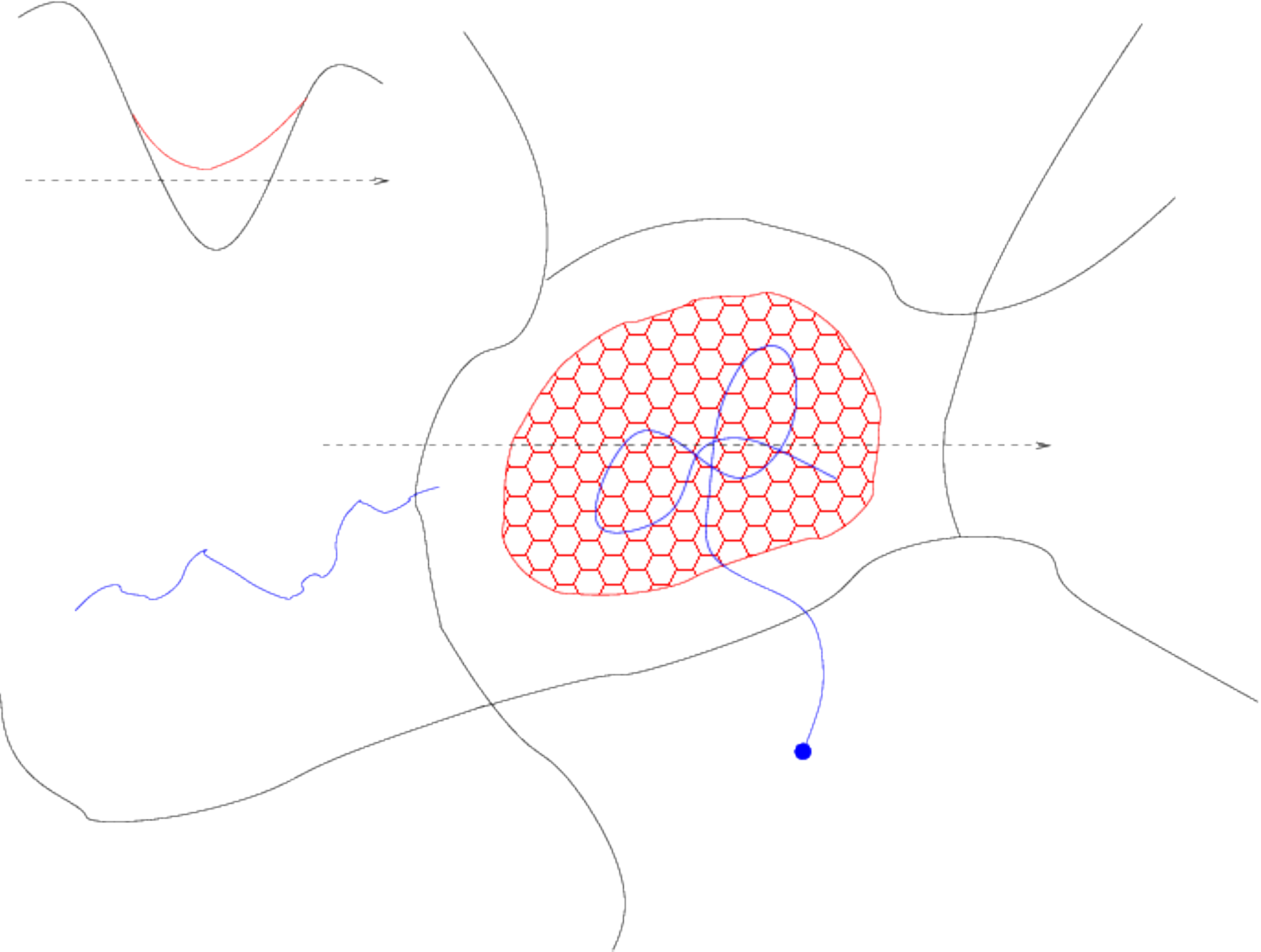}}
\caption{The hyperdynamics: the exit event is simulated on a biased
  potential in $S$. The honeycomb region represents the support of the
biasing potential $\delta V$.}\label{fig:Hyper}
\end{figure}

Let us again assume that we are given an initial condition $X_0 \in S$ such
that $X_0$ is distributed according to the QSD in $S$. In other words,
let us assume that we are at the end of the decorrelation step: the
reference walker stayed sufficiently long in $S$.

The principle
of the hyperdynamics algorithm is then to raise the potential inside the
state in order to accelerate the exit from $S$. The algorithm
thus requires a biasing potential $\delta V: S \to \R$, which
satisfies appropriate assumptions detailed below. The algorithm then
proceeds as follows:
\begin{itemize}
\item Equilibrate the dynamics on the biased potential $V+\delta V$,
  namely run the dynamics~\eqref{eq:overdamped_Lang} on the process
  $(X^{\delta V}_t)_{t \ge 0}$ over the biased
  potential conditionally to staying in the well, up to the time the
  random variable $X^{\delta V}_t$ has distribution close to the QSD
  $\nu^{\delta V}$ associated
  with the biased potential. This first step is a preparation step, which is pure
  overhead.  The end point $X^{\delta V}_t$ will be used as the
  initial condition for the next step.
\item Run the dynamics~\eqref{eq:overdamped_Lang} over the biased
  potential $V + \delta V$ up to the exit time $T^{\delta V}_S$ from
  the state $S$. The simulation clock is updated by adding the
  effective exit time $B T^{\delta V}_S$ where $B$ is the so-called
  boost factor defined by
\begin{equation}\label{eq:B}
B=\frac{1}{T^{\delta V}_S} \int_0^{T^{\delta V}_S} \exp(\beta \,
\delta V(X^{\delta V}_t)) \, dt.
\end{equation}
The exit point is then used as the starting point for a new
decorrelation step.
\end{itemize}
See
Figure~\ref{fig:Hyper} for a schematic illustration of the method.

Roughly speaking, the assumptions required on $\delta V$ in the
original paper~\cite{voter-97} are twofold:
\begin{itemize}
\item $\delta V$ is sufficiently small so that the exit event from the
  state $S$ still satisfies the standard assumptions used for kinetic
  Monte Carlo models and transition state theory.
\item $\delta V$ is zero on (a neighborhood) of the boundary $\partial
  S$.
\end{itemize}
The derivation of the method relies on explicit formulas for the laws of
the exit time and exit point, using the transition state theory.
The aim of the mathematical analysis presented below is to give a
rigorous set of assumptions to make this algorithm consistent.

The algorithm we present here is actually slightly different from the
way it is introduced in the original paper~\cite{voter-97}. Indeed, in
the original version, the local equilibration steps (decorrelation
step and equilibration step on the biased potential) are omitted: it is assumed that the states are sufficiently metastable
(for both the original potential and the biased potential)
so that these local equilibrations are immediate. It would be
interesting to check if the modifications we propose here improve the
accuracy of the method.

Let us now discuss the mathematical foundations of this technique, and
in particular, a way to understand the formula~\eqref{eq:B} for the
boost factor. We actually need to compare two exit events. The first
one is the exit event for the original process $X_t$ following the
dynamics~\eqref{eq:overdamped_Lang}, starting from the QSD $\nu$
associated with the state~$S$ and the dynamics with potential $V$. The
second one is the exit event for the process $X^{\delta V}_t$ following the
dynamics~\eqref{eq:overdamped_Lang} on the biased potential $V +
\delta V$, starting from the QSD $\nu^{\delta V}$
associated with the state~$S$ and the dynamics with potential $V+\delta
V$. Referring to Proposition~\ref{prop:QSD_3},  comparing the two exit
events amounts to understanding how the first eigenvalue $\lambda_1$
and the normal derivative of the first eigenvector $\partial_n u_1$ are modified when changing
the potential from $V$ to $V+\delta V$. Let us denote $\lambda_1(V)$
(resp. $\lambda_1(V+\delta V)$)
and $\partial_n u_1(V)$ (resp. $\partial_n u_1(V+\delta V)$) the first eigenvalue and the normal derivative
when considering the original potential $V$ (resp. the biased
potential $V+\delta V$).
In~\cite{lelievre-nier-13}, we prove the following.
\begin{theorem}\label{th:hyper}
Let us make the following assumptions on $V$. We assume there exists
an open set $S^{-}$ such that $\overline{S^-} \subset S$ and:
\begin{itemize}
\item {\em Regularity}: $V$ and $V|_{\partial S}$ are Morse
  functions.
\item {\em Localization in $S^-$ of the eigenvectors associated with
    small eigenvalues}:
\begin{enumerate}[(i)]
\item $|\nabla V|
  \neq 0$ in ${\overline S} \setminus S^-$\,;
\item   $\partial_n V > 0$ on $\partial  S^-$\,;
\item $\min_{\partial S} V \ge \min_{\partial S^-} V$\,; 
\item $\min_{\partial S^-} V - {\rm cvmax} > {\rm cvmax} - \min_{S^-} V$
where ${\rm cvmax}=\max\{V(x), x \text{ s.t. } |\nabla V(x)|=0 \}$\,.
\end{enumerate}
\item {\em Non degeneracy of exponentially small eigenvalues}: The critical values of $V$ in $S^-$ are all distinct and the
  differences $V(y) - V(x)$ are all distinct, where $x \in {\mathcal U}^{(0)}$ ranges over the local
  minima of $V|_{S^-}$ and $y \in {\mathcal U}^{(1)}$ ranges over the
  critical points of $V|_{S^-}$ with index $1$.
\end{itemize}
Let us also assume that the biasing potential $\delta V$ is such that
\begin{itemize}
\item $V + \delta V$ satisfies the same assumptions as the ones on $V$
  above ;
\item$\delta V = 0$ on $\overline S \setminus S^-$.
\end{itemize}
Then, there exists $c>0$ such that, in the limit $\beta \to \infty$,
\begin{align}
&
\frac{\lambda_{1}(V + \delta V)}{\lambda_{1}(V)}=\frac{\int_{S}
  e^{- \beta V}}{\int_{S} e^{-\beta (V + \delta V)}}
  (1+\mathcal{O}(e^{-\beta c}))\,, \label{eq:lambda_Hyp}\\
&\frac{\partial_{n}\left[u_{1}(V + \delta V)\right]\big|_{\partial S}}{ 
\left\|\partial_{n}\left[u_{1}(V + \delta V)\right]
\right\|_{L^{1}(\partial S)}}
= 
\frac{\partial_{n}\left[u_{1}(V)\right]\big|_{\partial S}}{ 
\|\partial_{n}\left[u_{1}(V)\right]\|_{L^{1}(\partial
 S)}} +
\mathcal{O}(e^{-\beta c})\quad \text{in}~ L^{1}(\partial S)\,. \label{eq:dnu_Hyp}
\end{align}
\end{theorem}
The proof is based on results from semi-classical analysis for
boundary Witten Laplacians.

Notice that the formula~\eqref{eq:lambda_Hyp} provides a justification
of the formula~\eqref{eq:B} for the boost factor. Indeed, by assuming
that $T^{\delta V}_S$ is sufficiently large, we have by an ergodic
property:
\begin{align*}
B&=\frac{1}{T^{\delta V}_S} \int_0^{T^{\delta V}_S} \exp(\beta \,
\delta V (X^{\delta V}_t)) \, dt.\\
&\simeq \frac{\int_S \exp(\beta \delta V) \exp(-\beta(V+\delta
  V))}{\int_S\exp(-\beta(V+\delta V))}\\
&= \frac{\int_S \exp(- \beta V) }{\int_S\exp(-\beta(V+\delta V))}\\
&\simeq \frac{\lambda_{1}(V + \delta V)}{\lambda_{1}(V)}.
\end{align*}
By multiplying the exit time on the biased potential $V+\delta V$ by
$\frac{\lambda_{1}(V + \delta V)}{\lambda_{1}(V)}$, we indeed obtain
(in law) the exit time on the original potential $V$.
Moreover, the estimate~\eqref{eq:dnu_Hyp} shows that (up to exponentially
small errors in the limit of small temperature), the first exit point
from $S$ for the biased potential has the same
distribution as the first exit point from $S$ for the original potential (see
Equation~\eqref{eq:QSD_3} in Proposition~\ref{prop:QSD_3}).

A practical aspect we do not discuss here at all is the
effective construction of the biasing potential $\delta V$. In the
original article~\cite{voter-97}, A.F. Voter proposes a technique based on the Hessian
$\nabla^2 V$. A well-known method in the context of
materials science is the bond-boost method introduced in~\cite{miron-fichthorn-03}.

Notice that, contrary to the Parallel Replica method, the
hyperdynamics is, at least for our mathematical analysis, limited
to energetic barriers (see assumptions {\em (iii)} and {\em (iv)} in
Theorem~\ref{th:hyper}). On the other hand, for very high energetic
barriers, hyperdynamics is in principle much more efficient: the Parallel Replica
method only divides the exit time by $N$ (the number of replicas),
while for deep wells, the boost factor $B$ is very large.

\section{The Temperature Accelerated Dynamics}\label{sec:TAD}

Let us finally introduce the Temperature Accelerated Dynamics (TAD),
see~\cite{sorensen-voter-00}. Let us assume again that we are at the
end of the decorrelation step: the reference walker stayed
sufficiently long in $S$. The principle of TAD is to increase the
temperature (namely increase $\beta^{-1}$
in~\eqref{eq:overdamped_Lang}) in order to accelerate the exit from
$S$. The algorithm consists in
\begin{itemize}
\item Simulating many exit events from $S$ at high temperature, starting
  from the QSD at high temperature,
\item Extrapolating the high temperature exit events to low
  temperature exit events using the Arrhenius law.
\end{itemize}
As for the hyperdynamics algorithm, in
the original paper~\cite{sorensen-voter-00}, no equilibration step is
used: it is assumed that the states are sufficiently metastable at
both high and low temperatures so that the convergence to the QSD is
immediate. Let us now describe more precisely how the extrapolation
procedure is made.

\begin{figure}[htbp]
\centerline{\includegraphics[width=7cm]{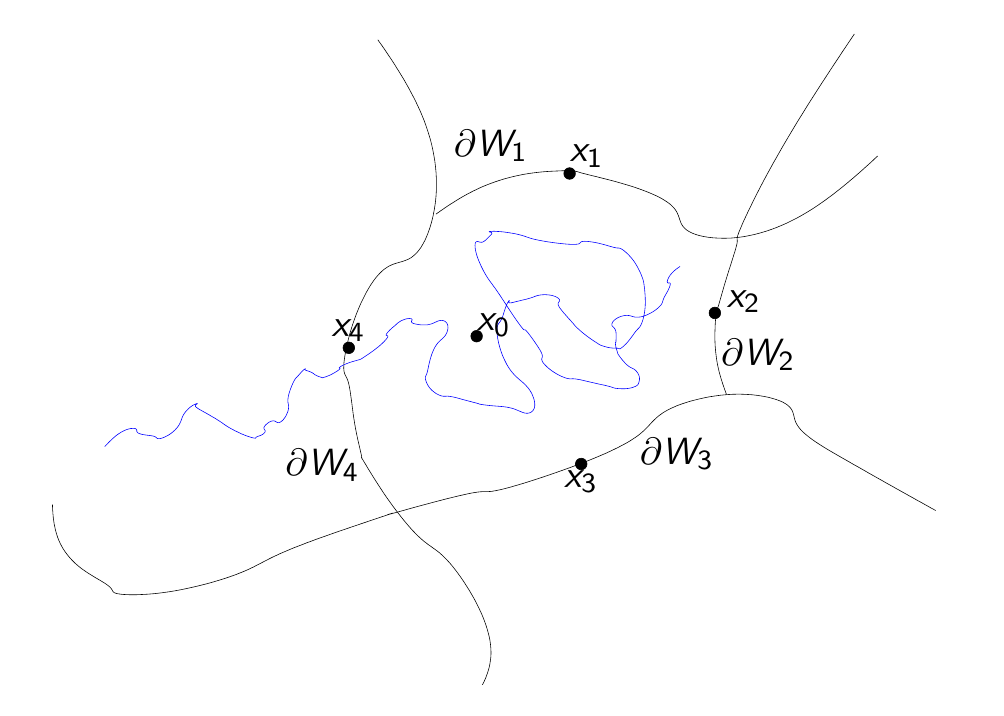}}
\caption{The Temperature Accelerated Dynamics: exit events are
  simulated at higher temperature and then extrapolated to the
  original smaller temperature.}\label{fig:TAD}
\end{figure}

Let us consider the exit event from $S$, at a given temperature. The
set $S$ is surrounded by $I$ neighboring states and let us denote by
$\partial S_i$ the common boundary with th $i$-th neighboring state,
$i \in \{1, \ldots, I\}$. The sets $\partial S_i$ thus form a
partition of the boundary $\partial S$. Let us introduce, for $i \in
\{1, \ldots, I\}$, the saddle point $x_i \in \partial S_i$ which is
the lowest in energy on $\partial S_i$: in the small temperature
regime, the paths leaving the state $S$ through $\partial S_i$ will
leave through a neighborhood of $x_i$ (this can be inferred from
results from the large deviation theory for example). Let us also
denote by $x_0$ the global minimum of $V$ on $S$. We refer to
Figure~\ref{fig:TAD} for a schematic representation of the
geometry. The interesting quantities to define the exit events are:
\begin{itemize}
\item The probability to exit through $\partial S_i$ which writes (see
  Proposition~\ref{prop:QSD_3}):
$$p(i)= \P(X_{T_S} \in \partial S_i)= - \frac{\displaystyle\int_{\partial S_i}
\partial_n u_1  \,
  d\sigma}{\displaystyle\beta \lambda_1 \int_{S} u_1(x)
  \, dx}.  $$
\item And the parameter of the exponential random variable $T_S$:
$$\lambda_1=1/\E(T_S).$$
\end{itemize}
Notice that one way to rewrite the exit event is to attach to each
exit regions $\partial S_i$ (or to each saddle point $x_i$) a rate
$$k(i)=\lambda_1 p(i)$$
and to consider $I$ independent exponential random variables $\tau_i$
with parameter $k(i)$. The exit event is then given by
\begin{itemize}
\item the exit time $\min(\tau_1, \ldots , \tau_I)\stackrel{\mathcal L}{=}T_S$
\item and the exit region $\arg \min(\tau_1, \ldots , \tau_I)$, since , for $i \in \{1,
  \ldots ,I\}$,
  $\P(\arg \min(\tau_1, \ldots , \tau_I)=i)=p(i)$.
\end{itemize}
This description of the exit event in terms of rates attached to neighboring
saddle points is exactly what is used for kinetic Monte Carlo
models~\cite{voter-05}. The TAD algorithm requires an approximation of the rate $k(i)$, namely the Arrhenius law:
\begin{equation}\label{eq:arrhenius}
k(i)=\lambda_1 p(i) \simeq \eta_i \exp(-\beta (V(x_i) - V(x_0)))
\end{equation}
where $\eta_i$ is independent of $\beta$. In the original
paper~\cite{sorensen-voter-00}, the Arrhenius law is justified using
the (harmonic) transition state theory.

Let us now go back to the TAD algorithm. Using the underlying kinetic Monte
Carlo model presented above, and using the Arrhenius
law~\eqref{eq:arrhenius}, one observes that:
\begin{equation}\label{eq:extrapolation}
\frac{k^{hi}(i)}{k^{lo}(i)}=\frac{\lambda^{hi}_1
  p^{hi}(i)}{\lambda^{lo}_1 p^{lo}(i)} \simeq
\exp(-(\beta^{hi}-\beta^{lo}) (V(x_i) - V(x_0)))
\end{equation}
where, with obvious notation, $\beta^{lo}$ denotes the inverse low
temperature, $\beta^{hi}$ the inverse high temperature, and the
superscripts $lo$ and $hi$ refer to the associated quantities
respectively at low and high temperature. The extrapolation
formula~\eqref{eq:extrapolation} is used in TAD in order to infer the exit event at
low temperature from the exit events observed at high temperature, by using the formula:
\begin{equation}\label{eq:TAD_justif}
(\tau_1^{lo},\ldots,\tau_I^{lo})\stackrel{\mathcal
  L}{=}(\Theta^1\tau_1^{hi},\ldots,\Theta^I\tau_I^{hi})
\end{equation}
where
$$\Theta^i=\frac{k^{hi}(i)}{k^{lo}(i)}\simeq\exp(-(\beta^{hi}-\beta^{lo}) (V(x_i) - V(x_0)))
$$
is a multiplicative factor constructed from the ratio of the
rates~\eqref{eq:extrapolation}. In the equality in law
in~\eqref{eq:TAD_justif} the random variables $\tau^{hi/lo}_i$ are, as
described above, exponential random variables with parameter
$k^{hi/lo}(i)$. To have analytical formula for the correction factors
$\Theta_i$ and make the algorithm practical, the Arrhenius is
assumed to be exact and one uses in practice
$\Theta^i=\exp(-(\beta^{hi}-\beta^{lo}) (V(x_i) - V(x_0)))$.

The TAD
algorithm thus consists in running the dynamics at high temperature,
observing the exit events through the saddle points on the boundary of
the state, and updating the exit time and exit region that would have
been observed at low temperature. More precisely, if exits through the
saddle points $\{s_1,\ldots,s_k\} \subset\{1, \ldots ,I\}$ have been
observed, one computes
$\min(\Theta^{s_1}\tau_{s_1}^{hi},\ldots,\Theta^{s_k}\tau_{s_k}^{hi})$ and
$\arg\min(\Theta^{s_1}\tau_{s_1}^{hi},\ldots,\Theta^{s_k}\tau_{s_k}^{hi})$
to get the exit time and the exit region at low temperature. 

The interest of TAD compared to a brute force saddle point search is
that it is not required to observe exits through all the saddle points
in order to obtain a statistically correct exit event. Indeed, a stopping criterium is introduced to
stop the calculations at high temperature when the extrapolation procedure
will not modify anymore the low temperature exit event
(namely will not modify $\min
(\Theta^{s_1}\tau_{s_1}^{hi},\ldots,\Theta^{s_k}\tau_{s_k}^{hi})$,
$\{s_1,\ldots,s_k\} \subset\{1, \ldots ,I\}$ being the saddle points
discovered up to the stopping time). This stopping
criterium requires to provide some a priori knowledge, typically a lower
bound on the barriers $V(x_i)-V(x_0)$ ($i \in \{1, \ldots, I\}$)
(there is also a variant using a lower
bound on the prefactors $\eta_i$ in~\eqref{eq:arrhenius}).
In some sense,
TAD can be seen as a clever saddle point search, with a rigorous way
to stop the searching procedure.

If the Arrhenius law~\eqref{eq:arrhenius} is exactly satisfied, one
can check that the TAD algorithm simulates an exit event which is
statistically exact. The mathematical question raised by this
algorithm is thus to estimate the difference between the ratio of the rates $\frac{\lambda^{hi}_1
  p^{hi}(i)}{\lambda^{lo}_1 p^{lo}(i)}$ and the estimate deduced from
the Arrhenius law $\exp(-(\beta^{hi}-\beta^{lo}) (V(x_i) - V(x_0)))$ (see the
extrapolation formula~\eqref{eq:extrapolation} above). In~\cite{aristoff-lelievre-14}, we
consider as a first step the case of a one dimensional potential,
where $S$ is a single well, and we prove that in the limit
$\beta^{hi},\beta^{lo} \to \infty$ with
$\beta^{lo}/\beta^{hi}$ fixed,
\begin{equation}\label{eq:error_TAD}
\frac{\lambda^{hi}p_i^{hi}}{\lambda^{lo}p_i^{lo}} = 
e^{-(\beta^{hi}-\beta^{lo})(V(x_i)-V(x_0))}\left(1 +
  O\left(\frac{1}{\beta^{hi}}- \frac{1}{\beta^{lo}}\right)\right).
\end{equation}
The extension of this result to a high dimensional setting is a work
under progress.

Notice that, compared to the hyperdynamics, TAD is based
on an additional assumption, namely the Arrhenius law. We expect that
this implies larger error for TAD than for the hyperdynamics: for the
hyperdynamics the error in~\eqref{eq:lambda_Hyp}--\eqref{eq:dnu_Hyp} is exponentially small in the
limit $\beta \to \infty$, while for TAD, the error
scales like $1/\beta$ in~\eqref{eq:error_TAD}. The interest of TAD 
compared to the hyperdynamics is that it does not require a
biasing potential, which may be complicated to build in some situation.

\section{Conclusion and discussion}\label{sec:conc}

We presented three algorithms which have been proposed by A.F. Voter
and co-workers in order
to efficiently generate the state-to-state dynamics associated with a
metastable stochastic process.

We proposed an analysis of these algorithms based on the notion of
quasi-stationary-distribution (QSD). As explained above, starting from
the QSD within a well, the exit event is Markovian since the exit time
is exponentially distributed and independent of the next visiting
state. From a theoretical viewpoint, the QSD thus seems to be an
interesting intermediate to relate Markovian dynamics in continuous
state space (such as the Langevin~\eqref{eq:Langevin} or the overdamped Langevin~\eqref{eq:overdamped_Lang} dynamics)
with Markovian dynamics in discrete state space (kinetic Monte Carlo
models or Markov state models). It also gives a natural definition of
a metastable region (see~\cite{lelievre-13}) as a region where the
stochastic process reaches local equilibrium (namely the QSD) before exiting.

Going from Parallel Replica to Hyperdynamics to TAD, the assumptions
required for the algorithm to be correct are more and more stringent. Indeed,
for Parallel Replica, no assumptions is required beyond the fact that
a good estimate of $\tau_{corr}$ is available (to assess the
convergence of the QSD). Hyperdynamics requires additional assumptions
on the potential: the metastability of the state should come from
energetic barriers (at least for our mathematical analysis to apply). Finally, the TAD algorithm requires in addition
the Arrhenius law to be satisfied. Likewise, going from Parallel
Replica to Hyperdynamics to TAD, the errors introduced by the algorithm
are expected to be larger and larger. On the other hand, the expected
computational gain with Parallel replica is typically smaller than for
the two other methods. In addition, comparing hyperdynamics with TAD, the interest
of TAD is that it does not require the construction of a biasing
potential, which is a difficult task in general.

One practical aspect we did not discuss above is the choice of the
partition of the configurational space into states (in other words, the choice of
the function ${\mathcal S}$). Let us focus on the Parallel Replica
algorithm for simplicity. We already mentioned that thanks to the
decorrelation step, the algorithm is consistent whatever the
choice of the partition. However, the efficiency of the algorithm
highly depends on the choice of the partition: the states should be
metastable regions, so that the stochastic process reaches the local
equilibrium (the QSD) before leaving the state. How to design a good
partition is a difficult question. For a system with high energy
barriers (this is often the case for application in materials science
for example), the original idea of A.F. Voter and co-workers is to
define the states as the basins of attraction of the local minima of
$\nabla V$ for the simple gradient dynamics $\dot{x}=-\nabla
V(x)$. For a system with more diffusive or entropic barriers (this is
typically the case for biological applications), one could think of
defining the states using relevant reaction coordinates (see for
example~\cite{kum-dickson-stuart-uberuaga-voter-04} where the states
  are defined in terms of the molecular topology of the molecule of
  interest). Notice that choosing a good partition also implies being
  able to estimate the correlation time within each state either from
  some {\em a priori} knowledge, or some on-the-fly estimates~\cite{binder-simpson-lelievre-14}.

Let us finally mention the mathematical questions raised by these
algorithms and that we would like to
investigate in the future.

Concerning the TAD algorithm, the analysis of the validity of the
Arrhenius law starting from the QSD has only been done for the moment
in a one-dimensional situation. The extension to a more general
setting is a work under progress.

We would like to stress that all these algorithms are used in practice
with the Langevin dynamics~\eqref{eq:Langevin}. The mathematical
results we presented above assumed that the dynamics was the
overdamped Langevin~\eqref{eq:overdamped_Lang}. Therefore, some works
have to be done to extend these results to the Langevin
dynamics. There are mainly two difficulties. First, the infinitesimal
generator associated with the overdamped Langevin is symmetric (in an
appropriate weighted $L^2$ space) while this is not the case for
Langevin, which implies that the study of the spectrum of the
infinitesimal generator for Langevin is more complicated. In addition, the fact that the domain of
interest is typically bounded in position but not in velocity implies
additional difficulties. We refer to
the recent work~\cite{nier-14} by F. Nier, which gives in
particular some conditions for the existence of an isolated smallest
eigenvalue (and therefore of the QSD) for the Fokker-Planck operator associated with  the Langevin dynamics.
 Second, while
there is an extensive literature on the spectrum of the infinitesimal
generator of the overdamped Langevin dynamics in the small temperature regime
(semi classical analysis for Schr\"odinger operators and Witten
laplacians), this is not the case for the Langevin dynamics. This
means that the analysis of the hyperdynamics or TAD will certainly be
more involved for Langevin than for overdamped Langevin.

\bibliography{biblio_HD,biblio_MD,ma_biblio}
\bibliographystyle{plain}

%

\end{document}